# Numerical Approximation In Real Domain Of Special Function Of Product Of A Variable And Its Double Exponential

**Narinder Kumar Wadhawan[1],**

[1]Civil Servant, Indian Administrative Service Retired,
House No. 563, Sector 2, Panchkula-134112, Haryana, India.
e-mails: narinderkw@gmail.com

*Abstract:* Purpose of writing this paper is to solve a transcendental function containing a product of a variable and its double exponential by a unique method of approximation. If the value of the said product is given, then its inverse function is approximated by use of linear expression in place of natural logarithm of a positive real quantity and, that transforms the function to a quadratic equation. Roots of the equation are, then used for solving the function. For precise approximation, the process is iterated a number of times and more the number of iterations, more precise will be the approximation. To prove truthfulness of the formulae derived, a number of examples are given in tabular form.
*Keywords*: Natural Logarithm; Quadratic Approximation; Real Number; Special Function.
*2010 Mathematics Subject Classifications:* 33E20, 33F05

## 1. Introduction

Equation $y\, e^{e^y} = X$, contains product of a variable $y$ and its double exponential $e^{e^y}$, where $X$ is a given real quantity except positive or negative infinity. A double exponential function is defined as a constant raised to the power of an exponential function. If constant is $a$ and $b^y$ is an exponential function, then double exponential function would be $a^{b^y}$, where $y$ is a variable. If $a = b = e$, where $e$ is an *Euler's number*, then the double exponential function takes the form $e^{e^y}$. Computation of such function is done by taking topmost part and, then moving downward. For example, $e^{e^y}$ is calculated, first by taking part $e^y$, then calculating it say $e^y = z$ and, thereafter, calculating $e^z$. Obviously, such function grows faster than normal exponential function. Such function is used in algorithm complexity [1], number theory [4],[5],[3], theoretical biology [6], physics [2] likewise. If the value of double exponential function is given, its inverse function can be found taking double logarithm that is logarithmic of logarithm of $x$.

If the function happens to be product of variable $y$ and double exponential function $e^{e^y}$, then value of $y$ can not be determined by simply resorting to double logarithm. This function is also not transformable to *Lambert W Function* $ye^y = x$. That means, *this function is unique in nature and is not transformable to a known special function.* It is submitted, in equation $y\, e^{e^y} = X$, in stead of taking $x$ in lower case, $X$ in upper case has been taken to avoid confusion at the time of establishing its relationship with *Lambert W Function* in section 5. In this paper, for given real value of $X$, this equation

$$ye^{e^y} = X, \qquad (1.1)$$

has been analysed for approximating the numerical value of $y$ and also for its other characteristics. On the analogy of *Lambert W Function*, this equation is given the nomenclature of $N$ function by writing it

$$N(X) = y \qquad (1.2)$$

and thus it is an inverse function of $ye^{e^y} = X$. It is submitted, in this paper, given $X$ is a real quantity. $X$ can also be a complex quantity of the form $a + i\,b$, where $a$ and $b$ are real quantities and $i = \sqrt{-1}$ but it is not the subject matter of this paper. Being double exponential, its properties obviously will not be the same as that of *Lambert W Function*. Considering $e^y = z$, the function can be written $(\ln z)e^z = X$ and that requires $z > 0$. This $N$ function can also be written $z\, \ln(\ln z) = X$, where $e^y = \ln z$ and that requires $z > 1$. That indicates, equations involving double exponential or logarithm and exponential or double logarithm can be solved provided the solution of Equation (1.1) is known. Here double logarithm is of a real quantity more than one. Notwithstanding the function as stated above, a general $N$ Function of the form
$$y^p e^{e^y} = X, \qquad (1.3)$$
has also been solved for $y$, for given real $p$ and $X$. *It is submitted, there are processes in physical and living sciences which involve i) a variable and double exponential of the variable or ii) a variable and double logarithm of the variable or iii) exponential of a variable and logarithm of the variable and such processes can be solved if solutions to Equations (1.1) and (1.3) are known. Denying the existence of such processes would infer as if all the processes that exist in the universe, are known and no such processes as stated above was found. Such a sweeping statement would tend to be more a fallacy than truthful. Solution of such unique equations underlines the importance of the present paper.* It will also be proved that the method of quadratic approximation used in this paper is self corrective meaning thereby that if an operator commits a mistake inadvertently, it will be self rectified at the cost of increased iterations. It will also be proved that initial assumption of value of $y$ has a wide range spectrum and choice of any value of $y$ out of the spectrum culminates in precise approximation same as for any other value of $y$. This characteristic property of the method is amusing and recreational. Icing in the cake is that the method used is unique and is quite different from traditional methods of Newton, Pade or any other approximation.

## 2. Plot Of The Function $N(X) = y$, When $X = y e^{e^y}$

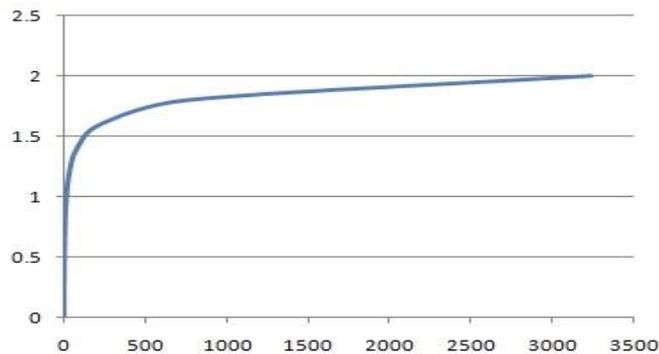

Figure 1 Plot of the function $N(X) = y$, when $X = y e^{e^y}$ and is more than 0

Referring to Figure 1, equation $y e^{e^y} = X$ is plotted by taking value of $X$ along x-axis from 0 to 3500 and taking corresponding value of $y$ along y-axis. It is observed, when value of $X$ is small, $y$ rises at a faster rate and when value of $X$ is large, the curve turns flat. Further, a single real value of $X$ corresponds to a single value of $y$ indicating, it is not a multivalued function.

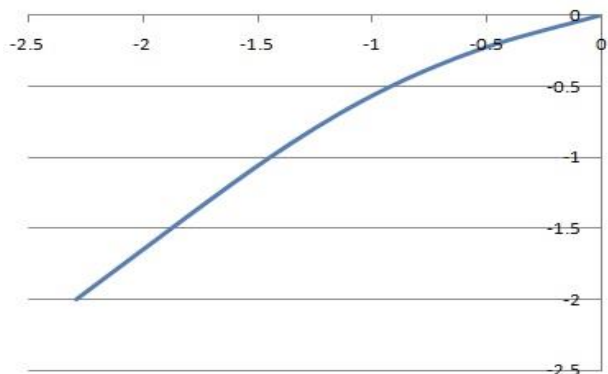

Figure 2 Plot of the function $N(-X) = -y$, when $-X = -y e^{e^{-y}}$ and $X > 0$

For negative values, of $X$, the curve takes the form as shown in Figure 2. At $X = 0$, value of $y$ is also zero. A given value of $X$ corresponds to one and only one value of $y$ in real domain.

## 2. Quadratic Approximation, Convergence, Rate Of Convergence And Associated Error

For numerical approximation of value of $y$, when value of $X$ is given, numerical approximation based on quadratic equation as explained in proceeding sections will be used.

### 2.1. Quadratic Equation Approximation

#### 2.1.1. First Method

For numerical approximation of value of $y$ in equation $ye^{e^y} = X$, it is written $z\,ln(ln\,z) = X$, where $e^y = ln\,z$ and value of $X$ is given. Value of $z$ is initially assumed $z_1$. How this assumption is made, is explained in section 3.1.1. Let the actual value of $z$ be $z_2 = z_1 + a_1$. To approximate the value of $a_1$, value of $z_1 + a_1$ is substituted for $z$ in equation $z\,ln(ln\,z) = X$ yielding

$$(z_1 + a_1)\,ln\{ln\,(z_1 + a_1)\} = X. \tag{2.1}$$

*It has been proved*

$$ln\left(\frac{t+1}{t}\right) \simeq \frac{2}{2t+1} \tag{2.2}$$

*provided $t \gg 1$ in a published paper* [7]. Equation (2.1), on applying Equation (2.2), transforms to

$$ln\,(ln\,z_1) + \frac{2a_1/\{ln\,z_1(2z_1 + a_1) + a_1\}}{2 + ln\,(ln\,z_1)(1 + ln\,z_1)} = X/(z_1 + a_1). \tag{2.3}$$

On simplification, this turns out to be a quadratic equation in $a_1$

$$a_1^2 - a_1 l_1 - m_1 = 0, \tag{2.4}$$

when

$$l_1 = -z_1 - \{2z_1 ln\,z_1\,ln\,(ln\,z_1) - X(1 + ln\,z_1)\}/\{2 + ln\,(ln\,z_1)(1 + ln\,z_1)\} \tag{2.5}$$

$$m_1 = -2z_1 ln\,z_1\{z_1 ln\,ln\,z_1 - X\}/\{2 + ln\,(ln\,z_1)(1 + ln\,z_1)\}. \tag{2.6}$$

Roots of this quadratic equations are

$$a_1 = (l_1/2) \pm (1/2)\sqrt{l_1^2 + 4m_1} \tag{2.7}$$

and root $a_1 = (l_1/2) + (1/2)\sqrt{l_1^2 + 4m_1}$ is considered for approximation as it yields

$$z_2 = z_1 + a_1 = z_1 + (l_1/2) + (1/2)\sqrt{l_1^2 + 4m_1} \tag{2.8}$$

and $z_n + a_n$ on repeated correction a positive quantity. For precise approximation, let $z = z_3 = z_2 + a_2$ where

$$a_2 = (l_2/2) + (1/2)\sqrt{l_2^2 + 4m_2}, \tag{2.9}$$

$$z_3 = z_2 + a_2 = z_2 + (l_2/2) + (1/2)\sqrt{l_2^2 + 4m_2}, \tag{2.10}$$

$$l_2 = -z_2 - \{2z_2 ln\,z_2\,ln\,n\,z_2) - X(1 + ln\,z_2)\}/\{2 + ln(ln\,z_2)(1 + ln\,z_2)\} \tag{2.11}$$

$$m_2 = -2z_2 ln\,z_2\{z_2 ln\,ln\,z_2 - X\}/\{2 + ln\,(ln\,z_2)(1 + ln\,z_2)\}. \tag{2.12}$$

The process is continued till precise approximation of $z = z_{n+1} = z_n + a_n$ is made, where $n$ is the number of times iteration is made and its value depends upon the extent of precision required. Finally,

$$a_n = (l_n/2) + (1/2)\sqrt{l_n^2 + 4m_n}, \tag{2.13}$$

$$z_{n+1} = z_n + a_n = z_n + (l_n/2) + (1/2)\sqrt{l_n^2 + 4m_n}, \tag{2.14}$$

$$l_n = -z_n - \{2z_n ln\,z_n\,ln\,(ln\,z_n) - X(1 + ln\,z_n)\}/\{2 + ln(ln\,z_n)(1 + ln\,z_n)\}, \tag{2.15}$$

$$m_n = -2z_n ln\,z_n\{z_n ln\,ln\,z_n - X\}/\{2 + ln\,ln\,z_n(1 + ln\,z_n)\}. \tag{2.16}$$

Value of $y$ then equals $ln\,\{ln\,(z_n + a_n)\}$.

#### 2.1.2. Second Method

In this method, it is considered, $e^y = z$ and equation $ye^{e^y} = X$ transforms to $(ln\,z)e^z = X$, which can be written as

$$z + ln(ln\,z) = ln\,X. \tag{2.17}$$

Let $z = z_1 + a_1$, then Equation (2.17) takes the form

$$z_1 + a_1 + \ln \ln (z_1 + a_1) = \ln X. \tag{2.18}$$

Equation (2.18) on application of Equation (2.2) and then its simplification, transforms to quadratic equation

$$a_1^2 - l_1 a_1 - m_1 = 0,$$

where

$$l_1 = -2(1 + z_1 \ln z_1)/(1 + \ln z_1) + (\ln X - \ln \ln z_1 - z_1) \tag{2.19}$$
$$m_1 = 2z_1 \ln z_1 (\ln X - \ln \ln z_1 - z_1)/(1 + \ln z_1). \tag{2.20}$$

and roots of this equation are $a_1 = (l_1/2) \pm (1/2)\sqrt{l_1^2 + 4m_1}$. Root

$$a_1 = (l_1/2) + (1/2)\sqrt{l_1^2 + 4m_1} \tag{2.21}$$

is considered for approximation for the reasons explained in section 2.1.1. What value $z_1$ should be assigned, has also been explained in forthcoming section 3.2.1. For precise approximation, the process is iterated $n$ times, by considering $z = z_2 + a_2$, where $z_2 = z_1 + a_1$, then $z = z_3 + a_3$, where $z_3 = z_2 + a_2$ and finally $z = z_n + a_n$, where $z_n = z_{n-1} + a_{n-1}$ and that will yield $y = \ln(z_n + a_n)$. Values of $a_2, l_2, m_2 \ldots a_n, l_n, m_n$ are given by equations

$$a_2 = (l_2/2) + (1/2)\sqrt{l_1^2 + 4m_2}, \tag{2.22}$$
$$l_2 = -2(1 + z_2 \ln z_2)/(1 + \ln z_2) + (\ln X - \ln \ln z_2 - z_2), \tag{2.23}$$
$$m_2 = 2z_2 \ln z_2 (\ln X - \ln \ln z_2 - z_2)/(1 + \ln z_2). \tag{2.24}$$

...

$$a_n = (l_n/2) + (1/2)\sqrt{l_n^2 + 4m_n}, \tag{2.25}$$
$$l_n = -2(1 + z_n \ln z_n)/(1 + \ln z_n) + (\ln X - \ln \ln z_n - z_n), \tag{2.26}$$
$$m_n = 2z_n \ln z_n (\ln X - \ln \ln z_n - z_n)/(1 + \ln z_n). \tag{2.27}$$

### 2.1.3. Third Method

In this method, $-ye^{e^{-y}} = -X$ is written $-y = \ln \ln(X/y)$. Let $y = y_1 + a_1$, then this equation takes the form

$$-(y_1 + a_1) = \ln \ln\{X/(y_1 + a_1)\} \tag{2.28}$$

Equation (2.28) on application of Equation (2.2) and its simplification, transforms to quadratic equation

$$a_1^2 - l_1 a_1 - m_1 = 0,$$

where

$$l_1 = -\{y_1 + \ln \ln(X/y_1)\} + 2\{y_1 \ln (X/y_1) - 1\}/\{1 - \ln (X/y_1)\}, \tag{2.29}$$
$$m_1 = 2y_1 \ln (X/y_1)\{y_1 + \ln \ln (X/y_1)\}/\{1 - \ln (X/y_1)\} \tag{2.30}$$

and roots of this equation are $a_1 = (l_1/2) \pm (1/2)\sqrt{l_1^2 + 4m_1}$. Root

$$a_1 = (l_1/2) + (1/2)\sqrt{l_1^2 + 4m_1} \tag{2.31}$$

is considered for approximation for the reasons explained in section 2.1.1. What value $y_1$ should be assigned, has been explained in section 3.3.1. For precise approximation, the process is iterated $n$ times, by considering $y = y_2 + a_2$, where $y_2 = y_1 + a_1$, then $y = y_3 + a_3$, where $y_3 = y_2 + a_2$ and finally $y = y_n + a_n$, where $y_n = y_{n-1} + a_{n-1}$. Values of $a_2, l_2, m_2 \ldots a_n, l_n, m_n$ will be given by equations

$$a_2 = (l_2/2) + (1/2)\sqrt{l_1^2 + 4m_2}, \tag{2.32}$$
$$l_2 = -\{y_2 + \ln \ln(X/y_2)\} + 2\{y_2 \ln (X/y_2) - 1\}/\{1 - \ln (X/y_2)\}, \tag{2.33}$$
$$m_2 = 2y_2 \ln (X/y_2)\{y_2 + \ln \ln (X/y_2)\}/\{1 - \ln (X/y_2)\}, \tag{2.34}$$

...

$$a_n = (l_n/2) + (1/2)\sqrt{l_n^2 + 4m_n}, \tag{2.35}$$
$$l_n = -\{y_n + \ln \ln(X/y_n)\} + 2\{y_n \ln (X/y_n) - 1\}/\{1 - \ln (X/y_n)\}, \tag{2.36}$$
$$m_2 = 2y_n \ln (X/y_n)\{y_2 + \ln \ln (X/y_n)\}/\{1 - \ln (X/y_n)\}. \tag{237}$$

Value of $y$ thus equals $y_n + a_n$.

### 2.2. Convergence

### 2.2.1. Convergence of z approximated from equation $z \ln(\ln z) = X$

For examining the convergence of $z$ in equation $z \ln(\ln z) = X$ obtained from $y e^{e^y} = X$, where $e^y = \ln z$, this equation $z \ln(\ln z) = X$ is transformed to quadratic Equation (2.4) and value of its root $a_1$ is given by Equation (2.7). Kindly refer to section 2.1.1. Value of $z = z_2$ obtained from Equation 2.8 will make *approximation of $z \ln(\ln z)$ with X more precise than value of $z = z_1$ on account of quadratic equation approximation based on Equation 2.2* [7]. Hence magnitude of $|z_2 \ln(\ln z_2) - X| < |z_1 \ln(\ln z_1) - X|$. Assuming $z_2 > e$, value of $l_2$, according to Equation (2.11), will have negative value, whereas value of $m_2$ according to Equation (2.12) will be positive or negative depending upon whether $z_2 \ln(\ln z_2)$ is less than $X$ or greater than $X$. Let $l_2 = -L_2$, then value of $a_2$ on simplification is given by relation $a_2 = -m_2/L_2 = m_2/l_2$ since magnitude of $m_2 \ll l_2$. Quantity $|z_2 \ln(\ln z_2) - X|$, being small, on account of first correction, therefore on simplification,

$$|a_2| = |\ln z_2 \{z_2 \ln(\ln z_2) - X\}/\{1 + \ln z_2 \ln(\ln z_2)\}| \tag{2.38}$$

Value of $z$ will, then be assumed equal to $z_3 = z_2 + a_2$ and that will make approximation of $z \ln(\ln z)$ with $X$ more precise than value of $z = z_2 = z_1 + a_1$ on account of quadratic equation approximation based on Equation 2.2 [7]. Therefore, $|z_3 \ln(\ln z_3) - X| < |z_2 \ln(\ln z_2) - X|$ and since

$$|a_3| = |\ln z_3 \{z_3 \ln(\ln z_3) - X\}/\{1 + \ln z_3 \ln(\ln z_3)\}|,$$

therefore, $|a_3| < |a_2|$. Proceeding in this manner, it is proved,

$$|a_n| < |a_{n-1}| < |a_{n-2}| < \cdots < |a_2| < |a_1|$$

Hence $z$ converges to value $(z_{n+1})$, where $z_{n+1} = z_1 + a_1 + a_2 + a_3 + \cdots + a_n$ and accordingly, $y$ converges to value $\ln \ln (z_{n+1})$.

Rate of convergence is a measure that depicts how fast the sequence converges [8]. From Equation (2.38), value of $a_k$ is written

$$|a_k| = |\ln z_k \{z_k \ln(\ln z_k) - X\}/\{1 + \ln z_k \ln(\ln z_k)\}|,$$

where $k$ is a positive integer and as explained above, series $a_k$ is convergent, therefore,

limit $k$ tending to infinity, $a_{k+1}$ say $a_L$ tends to zero,.

Thus limit $k$ tending to infinity,

$$|(a_{k+1} - a_L)|/|(a_k - a_L)| = |(a_{k+1} - 0)|/|(a_k - 0)| = |a_{k+1}|/|a_k|$$

As proved above, $|a_{k+1}| < |a_k|$, $a_{k+1}$ tends to zero, when $k$ tends to infinity and $a_k$ has appreciably small value approaching to zero, therefore, $|a_{k+1}|/|a_k| = 0$. That means, when k tends to infinity,

$$|(a_{k+1} - a_L)|/|(a_k - a_L)| = 0.$$

*That proves sequence $a_k$ converges super linearly to zero [8]*

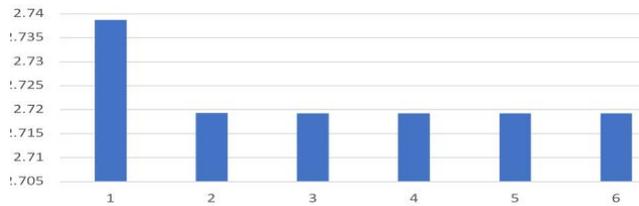

Figure 3 Convergence of $z$ in $z \ln(\ln z) = X$. when $z = z_1 + a_1 + a_2 + \cdots + a_6$, $X = 10^{-3}$ and $z_1 = 2$

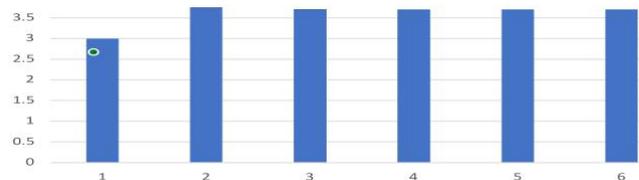

Figure 4 Convergence of $z$ in $z \ln(\ln z) = X$. when $z = z_1 + a_1 + a_2 + \cdots + a_6$, $X = 1$ and $z_1 = 2$

Figures 3 and 4 are drawn taking values of $z_2, z_3, ..., z_6$ from Table 1. It is observed from the Figure that length of last vertical bar $z_7$ approximates length of preceding bar $z_6$.

*2.2.2. Convergence of z approximated from equation $e^z \ln z = X$*

For examining the convergence of $y$ in equation $ye^{e^y} = X$, it is written as $e^z \ln z = X$, where $e^y = z$, this equation $e^z \ln z = X$, is further written as $z + ln(ln\ z) = ln\ X$ and convergence of z from this equation will be examined. Referring to section 2.1.2, this equation is transformed to quadratic equation $a_1^2 - l_1 a_1 - m_1 = 0$ and value of its root $a_1$ is given by Equation (2.21). Value of $z = z_2$ obtained by Equation 2.22 will make approximation of $z + \ln(ln\ z)$ with $ln\ X$ more precise than value of $z = z_1$ *on account of quadratic equation approximation based on Equation 2.2 [7]*. Hence *magnitude of* $|\{z_2 + ln(ln\ z_2) - ln\ X\}| < |\{z_1 + ln(ln\ z_1) - ln\ X\}|$. Assuming $z > 1$ and also $z_2 > 1$, value of $l_2$ according to Equation 2.23, will have negative value, since $-2(1 + z_2 ln\ z_2)/(1 + ln\ z_2)$ will be less than $-2$ and $|\{z_2 + \ln(ln\ z_2) - ln\ X\}|$ will be approaching zero, whereas value of $m_2$ according to Equation (2.24) will be positive or negative depending upon whether $z_2 \ln(ln\ z_2) + z_2 < X$ or $z_2 \ln(ln\ z_2) + z_2 > X$. Let $l_2 = -L_2$ then value of $a_2$ on simplification is given by relation $a_2 = -m_2/L_2 = m_2/l_2$, since magnitude of $m_2 \ll l_2$. Also quantity $|(ln\ X - ln\ ln\ z_2 - z_2)|$, being small on account of first correction, therefore, on simplification,

$$|a_2| = |z_2 ln\ z_2 \{ln\ X - z_2 - ln(ln\ z_2)\}/(1 + z_2 ln\ z_2)|$$

Now value of $z = z_3 = z_2 + a_2$ will make approximation of $\{z_2 + ln(ln\ z_2)\}$ with $ln\ X$ more precise than value of $z = z_2 = z_1 + a_1$ *on account of quadratic equation approximation based on Equation 2.2 [7]*. Therefore, $|\{z_3 + ln(ln\ z_3) - ln\ X\}| < |\{z_2 + ln(ln\ z_2) - ln\ X\}|$ and since

$$|a_3| = |z_3 ln\ z_3 \{ln\ X - z_3 - ln(ln\ z_3)\}/(1 + z_3 ln\ z_3)|$$

therefore, $|a_3| < |a_2|$. Proceeding in this way, it is proved,

$$|a_n| < |a_{n-1}| < |a_{n-2}| < \cdots < |a_2| < |a_1|$$

Hence z converges to $z_{n+1}$, therefore, y converges to value $ln(z_{n+1})$, where

$$z_{n+1} = z_1 + a_1 + a_2 + a_3 + \cdots + a_n.$$

For deciding rate of convergence, value of $a_k$ is written

$$|a_k| = |z_k ln\ z_k \{ln\ X - z_k - ln(ln\ z_k)\}/(1 + z_k ln\ z_k)|,$$

where $k$ is a positive integer [8] and as explained above, series $a_k$ is convergent, therefore,

limit $k$ tending to infinity, $a_{k+1}$ say $a_L$ tends to zero,.

Thus limit $k$ tending to infinity,

$$|(a_{k+1} - a_L)|/|(a_k - a_L)| = |(a_{k+1} - 0)|/|(a_k - 0)| = |a_{k+1}|/|a_k|$$

As proved above, $|a_{k+1}| < |a_k|$, $a_{k+1}$ tends to zero, when $k$ tends to infinity and $a_k$ has appreciably small value approaching to zero, therefore, $|a_{k+1}|/|a_k| = 0$. That means, when k tends to infinity,

$$|(a_{k+1} - a_L)|/|(a_k - a_L)| = 0.$$

*That proves sequence $a_k$ converges super linearly to zero [8].*

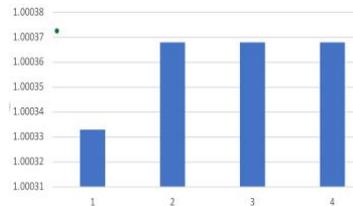

Figure 5 Convergence of z in $z + ln(ln\ z) = ln\ X$. when $z = z_1 + a_1 + a_2 + a_3 + a_4$, $X = 10^{-3}$ and $z_1 = 1.001$

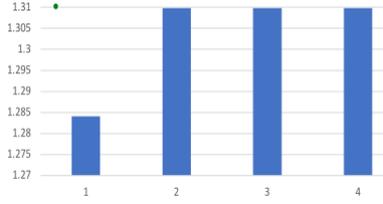

Figure 6 Convergence of $z$ in $z + ln(\ln z) = \ln X$. when $z = z_1 + a_1 + a_2 + a_3 + a_4$, $X = 1$ and $z_1 = 2$

Figures 5 and 6 are drawn taking values of $z_2, z_3, z_4, z_5$ from Table 2. It is observed from the Figures that length of last vertical bar $z_5$ approximates length of preceding bar $z_4$.

### 2.2.3. Convergence of $y$ in equation $-ye^{e^{-y}} = -X$

For examining the convergence of $y$ in equation $-ye^{e^{-y}} = -X$, the equation is written as $y = ln\ ln(X/y)$. Referring to section 2.1.3 and assuming $y = y_1 + a_1$, this equation is transformed to quadratic equation $a_1^2 - l_1 a_1 - m_1 = 0$ and value of its root $a_1$ is given by Equation (2.31). Value of $y = y_2 = y_1 + a_1$ will make approximation of $-y_2$ with $ln\ ln(X/y_2)$ more precise than value of $y = y_1$ *on account of quadratic equation approximation based on Equation 2.2* [7]. Hence magnitude of $|y_2 + ln \ln (X/y_2)| < |y_1 + ln \ln (X/y_1)|$. Assuming value of $y_1$ such that $X/y_1 > 1$, then value of $y_2$ obtained by quadratic equation approximation will yield value of $l_2$ as given by Equation (2.33) a negative quantity and value of $m_2$ given by Equation (2.34) will be a positive or negative depending upon whether $\{y_2 + ln \ln (X/y_2)\}$ is greater than or less than zero. Since value of $y = y_2$ will make approximation of $-y_2$ with $ln\ ln(X/y_2)$ more precisely, therefore, value of $\{y_2 + ln\ ln(X/y_2)\}$ will be small approaching to zero.. Value of $|a_2|$ will equal $|m_2/l_2|$ as has been explained in sections 2.2.1 and 2.2.2. Substitution of values of $l_2$ and $m_2$ as given by Equations (2.33) and (2.34) in equation $|a_2| = |m_2/l_2|$ and then its simplification yield,

$$|a_2| = |y_2 ln\ (X/y_2)\{y_2 + ln \ln (X/y_2)\}/\{y_2 \ln (X/y_2) - 1\}|.$$

Since $|y_2 + ln \ln (X/y_2)| < |y_1 + ln \ln (X/y_1)|$, therefore, $|a_2| < |a_1|$. Similarly, it can be proved that $|a_3| < |a_2|$, $|a_4| < |a_3|$ so on and $|a_n| < |a_{n-1}|$ or

$$|a_n| < |a_{n-1}| < |a_{n-2}| < \cdots < |a_2| < |a_1|.$$

Hence $y$ converges to value $y_{n+1}$, where $y_{n+1} = y_1 + a_1 + a_2 + a_3 + \cdots + a_n$.

For deciding rate of convergence, value of $a_k$ is written

$$|a_k| = |y_k ln\ (X/y_k)\{y_k + ln \ln (X/y_k)\}/\{y_k \ln (X/y_k) - 1\}|,$$

where $k$ is a positive integer and as explained above, series $a_k$ is convergent, therefore,

limit $k$ tending to infinity, $a_{k+1}$ say $a_L$ tends to zero,

Thus limit $k$ tending to infinity,

$$|(a_{k+1} - a_L)|/|(a_k - a_L)| = |(a_{k+1} - 0)|/|(a_k - 0)| = |a_{k+1}|/|a_k|$$

As proved above, $|a_{k+1}| < |a_k|$, $a_{k+1}$ tends to zero, when $k$ tends to infinity and $a_k$ has appreciably small value approaching to zero, therefore, $|a_{k+1}|/|a_k| = 0$. That means, when k tends to infinity,

$$|(a_{k+1} - a_L)|/|(a_k - a_L)| = 0.$$

*That proves sequence $a_k$ converges super linearly to zero.*

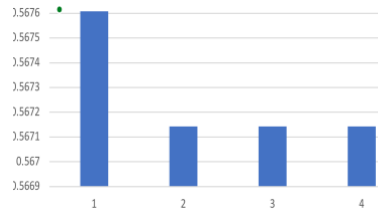

Figure 7 Convergence of $y$ in $y = ln\ ln(X/y)$, when $y = y_1 + a_1 + a_2 + a_3 + a_4$, $X = 1$ and $y_1 = 0.5$

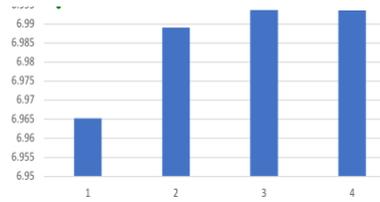

Figure 8 Convergence of $y$ in $y = \ln \ln(X/y)$, when $y = y_1 + a_1 + a_2 + a_3 + a_4$, $X = 7$ and $y_1 = 6.999$

Figures 7 and 8 are drawn taking values of $y_2, y_3, y_4, y_5$ from Table 3. It is observed from the Figures that length of last vertical bar $y_5$ approximates length of preceding bar $y_4$.

### 2.3. Error

It has been observed from the Tables that will follow that error associated with Quadratic Equation Approximation vanishes, when precision up to nine decimal points is required and iteration is made six times. For still better precision, further iterations and high precision calculator will be needed. For this paper, Desmos Scientific calculator precise up to nine decimal points has been used. In strictest sense, percentage error $E$ will be given by equation

$$E = 100(Y - y_{n+1})/y = 100(a_{n+1} + a_{n+2} + a_{n+3} + \cdots + a_{infinity})/Y,$$

where $Y$ is actual value of $y$ that satisfies equation $ye^{e^y} = X$. It is evident, higher the value of $n$, less would be the percentage error.

### 3. Numerical Approximation Of $y$ In Equation $ye^{e^y} = X$

Two methods to numerically approximate the value of $y$ in equation $ye^{e^y} = X$, are given hereinafter in sections 3.1 and 3.2, when $\infty > X > 0$ and third method is given in section 3.3 for equation $-ye^{e^{-y}} = -X$, when $\infty > X > 0$.

### 3.1. Numerical Approximation Of $y$, Where $y = \ln(\ln z)$, In Equation $z \ln(\ln z) = X$

For explaining numerical approximation of $y$, reference is made to section 2.1.1 whereby equation $ye^{e^y} = X$ is written as $z \ln(\ln z) = X$, where $y = \ln(\ln z)$ and $y$ is a positive real quantity. This equation is transformed to quadratic Equation (2.4) and its root given by Equation (2.7) is considered for approximating value of $z = z_1 + a_1$, where $z_1 > 1$. Further correction is applied by way of $a_2$ which is given by Equation (2.9) and $z = z_3 = z_2 + a_2$. In this way, finally $a_n$ is given by Equation (2.13) and $z = z_{n+1} = z_n + a_n$. Value of $y$, then equals $\ln(\ln z_{n+1})$. To avoid repetition, the details are not given in this section but these can be had by referring to section 2.1.1.

#### 3.1.1. Assignment Of Value To $z_1$

Since $y > 0$, it implies value of $z_1$ will be more than Euler number $e$. This method of approximation is self correcting in nature as is explained in forthcoming sections, even assumption of value less than Euler's number $e$ but more than 1 will yield correct result provided values of $(z_1 + a_1), (z_2 + a_2), (z_3 + a_2), \ldots, (z_n + a_n)$ are all more than one. Kindly refer to Table 1. Here assumption of value of $z_1 = 2$ could solve equation $z \ln(\ln z) = X$, when $X$ varies from $10^{-3}$ to $10^5$. It is also explicit from Table 4 that $z_1$ can be assigned any value from 2 to $10^{10}$ for solving equation $z \ln(\ln z) = 10^3$. In fact $z_1$ can be assigned any value from 2 to $10^p$ provided calculator is precise enough to detect $X$, when dealing with calculations of number as high as $10^p$, where $p$ is a positive integer such that $10^p \gg X$. In case of limited precision of Desmos Scientific calculator being used for this paper, for $X = 10^3$, maximum value of $z_1$ that can be assigned was found $10^{10}$.

On the basis of the formulae derived in section 2.1.1, the approximation of value $y$ is given in Table 1. No error was found in approximated value of $y$ up to nine decimal points, when correction was applied six times or less. This proves Lemma 3.1.

Table 1 Displaying correctness of quadratic equation approximation of z in $z \ln(\ln z) = X$

| $X$ | $z_1$ | $z_2$ | $z_3$ | $z_4$ | $z_5$ | $z_6$ | $z_7$ | $y = \ln(\ln z_8)$ | Actual |
|---|---|---|---|---|---|---|---|---|---|
| $10^{-3}$ | 2 | 2.73872599 | 2.719282203 | 2.71928183 | 2.719281829 | 2.7192818437 | No need | $3.67744156 \times 10^{-4}$ | $3.67744156 \times 10^{-4}$ |
| 1 | 2 | 2.998335266 | 3.760534672 | 3.710149812 | 3.704249825 | 3.705724821 | 3.705357822 | .2698741376 | .2698741376 |
| 10 | 2 | 14.7035688 | 11.28680512 | 11.29308446 | No need | No need | No need | 0.885497672 | 0.885497672 |
| $10^5$ | 2 | 122744.0354 | 42122.72914 | 42270.20703 | No need | No need | No need | 2.36573244 | 2.36573244 |

**Lemma 3.1**: *Numerical value of $y$ in equation $ye^{e^y} = X$, when $\infty > X > 0$ can be approximated by quadratic equation using relation $y = \ln \ln (z_n + a_n)$, where $z_1$ is assumed 2, $n$ pertains to number of times the correction is applied depending upon the extent of accuracy required and*

$$z_1 + a_1 = z_1 + (1/2)\left(l_1 + \sqrt{l_1^2 + 4m_1}\right),$$
$$l_1 = -z_1 - \{2z_1 \ln z_1 \ln (\ln z_1) - X(1 + \ln z_1)\}/\{2 + \ln (\ln z_1)(1 + \ln z_1)\}$$
$$m_1 = -2z_1 \ln z_1 \{z_1 \ln \ln z_1 - X\}/\{2 + \ln (\ln z_1)(1 + \ln z_1)\}.$$
$$z_2 + a_2 = z_2 + (1/2)\left(l_2 + \sqrt{l_2^2 + 4m_2}\right),$$
$$l_2 = -z_2 - \{2z_2 \ln z_2 \ln n z_2) - X(1 + \ln z_2)\}/\{2 + \ln(\ln z_2)(1 + \ln z_2)\}$$
$$m_2 = -2z_2 \ln z_2 \{z_2 \ln \ln z_2 - X\}/\{2 + \ln (\ln z_2)(1 + \ln z_2)\}.$$
$$\ldots$$
$$z_n + a_n = z_n + (1/2)\left(l_2 + \sqrt{l_n^2 + 4m_n}\right),$$
$$l_n = -z_n - \{2z_n \ln z_n \ln (\ln z_n) - X(1 + \ln z_n)\}/\{2 + \ln(\ln z_n)(1 + \ln z_n)\},$$
$$m_n = -2z_n \ln z_n \{z_n \ln \ln z_n - X\}/\{2 + \ln \ln z_n (1 + \ln z_n)\}.$$

### 3.2. Numerical Approximation Of $y$, Where $y = \ln z$, In Equation $z + \ln(\ln z) = \ln X$.

For explaining numerical approximation of $y$, reference is made to section 2.1.2, whereby equation $ye^{e^y} = X$, *when* $\infty > X > 0$, is written as $z + \ln(\ln z) = \ln X$, where $y = \ln z$ is a positive real quantity. This equation is transformed to quadratic equation $a_1^2 - l_1 a_1 - m_1 = 0$ and its root given by Equation (2.21) is considered for approximating value of $z = z_1 + a_1$, where $z_1 > 1$. Further correction is applied by way of $a_2$ which is given by Equation (2.22) and $z = z_3 = z_2 + a_2$. In this way, finally $a_n$ is given by Equation (2.25) and $z = z_{n+1} = z_n + a_n$. Value of $y$ then equals $\ln z$. To avoid repetition, the details are not given in this section but these can be had by referring to section 2.1.2.

### 3.2.1. Assignment Of Value To $z_1$

Since $y > 0$, it implies $z_1$ should be more than one. This method of approximation is self corrective in nature as is explained in forthcoming sections. Assumption of value more than 1 will yield correct result provided values of $(z_1 + a_1), (z_2 + a_2), (z_3 + a_2), \ldots, (z_n + a_n)$ are all more than one. Kindly refer to Table 2. Here assumption of value of $z_1$ equal to 2 could solve equation $z + \ln(\ln z) = X$, when $X$ varies from 1 to $10^{20}$. But when $X = 10^{-3}$, $z_1$ was assumed equal to 1.001. It is also explicit from Table 5 that $z_1$ can be assigned any value from 2 to $10^8$ for solving equation $z + \ln(\ln z) = \ln(10^2)$. In fact $z_1$ can be assigned any value from 2 to $10^p$ provided calculator is precise enough to detect value of $\ln (X)$ when dealing with calculations of number as high as $10^p$, where $p$ is a positive integer such that $10^p \gg \ln(X)$. In case of limited precision of Desmos Scientific calculator, being used for this paper, for $X = 10^2$, maximum value of $z_1$ that can be assigned was found $10^8$.

On the basis of the formulae derived in section 2.1.2, the approximation of value $y$ is given in Table 2. No error was found in approximated value of $y$ up to nine decimal points, when correction was applied four times or less. This proves Lemma 3.2.

Table 2 Displaying correctness of quadratic equation approximation of z in $z + \ln(\ln z) = X$

| $X$ | $z_1$ | $z_2 = z_1 + a_1$ | $z_3 = z_2 + a_2$ | $z_4 = z_3 + a_3$ | $z_5 = z_4 + a_4$ | $y = \ln z_5$ | Actual val. |
|---|---|---|---|---|---|---|---|
| $10^{-3}$ | 1.001 | 1.000333296 | 1.000367848 | 1.000367811 | 1.000367812 | $3.67744374 \times 10^{-4}$ | $3.67744156 \times 10^{-4}$ |
| 1 | 2 | 1.284087829 | 1.309809887 | 1.309799586 | 1.3097995858 | 0.2698741376 | .2698741376 |
| $10^2$ | 2 | 4.283737107 | 4.237733285 | 4.237733378 | No need | 1.444028546 | 1.44402854 |
| $10^5$ | 2 | 10.88208259 | 10.6518374 | 10.65183779 | No need | 2.36573244 | 2.3657324 |
| $10^{10}$ | 2 | 21.89883456 | 21.89883476 | No need | No need | 3.086433428 | 3.0864334272 |
| $10^{20}$ | 2 | 45.2800161 | 44.7166103 | 44.71661001 | No need | 3.800345021 | 3.800345021 |

**Lemma 3.2**: *Numerical value of $y$ in equation $ye^{e^y} = X$, when $\infty > X > 0$ can be approximated by quadratic equation*
*using relation $y = \ln(z_n + a_n)$, where $z_1$ is assumed 2, for $X \geq 1$ and $2 > z_1 > 1$ for $X \leq 1$, $n$ pertains to number of times the correction is applied depending upon the extent of accuracy required and*

$$z_1 + a_1 = z_1 + (1/2)\left(l_1 + \sqrt{l_1^2 + 4m_1}\right),$$
$$l_1 = -2(1 + z_1 \ln z_1)/(1 + \ln z_1) + (\ln X - \ln \ln z_1 - z_1)$$
$$m_1 = 2z_1 \ln z_1 (\ln X - \ln \ln z_1 - z_1)/(1 + \ln z_1).$$
$$z_2 + a_2 = z_2 + (1/2)\left(l_2 + \sqrt{l_2^2 + 4m_2}\right),$$
$$l_1 = -2(1 + z_1 \ln z_1)/(1 + \ln z_1) + (\ln X - \ln \ln z_1 - z_1)$$
$$m_1 = 2z_1 \ln z_1 (\ln X - \ln \ln z_1 - z_1)/(1 + \ln z_1).$$
$$\ldots$$
$$z_n + a_n = z_n + (1/2)\left(l_2 + \sqrt{l_n^2 + 4m_n}\right),$$
$$l_n = -2(1 + z_n \ln z_n)/(1 + \ln z_n) + (\ln X - \ln \ln z_n - z_n),$$
$$m_n = 2z_n \ln z_n (\ln X - \ln \ln z_n - z_n)/(1 + \ln z_n).$$

### 3.3. Numerical Approximation Of $y$, In Equation $-ye^{e^{-y}} = -X$

For explaining numerical approximation of $y$, reference is made to section 2.1.3, whereby equation $-ye^{e^y} = -X$, when $\infty > X > 0$, is written as $-y + \ln(\ln y) = \ln(\ln X)$ where, $y$ is a positive real quantity. This equation is transformed to quadratic equation $a_1^2 - l_1 a_1 - m_1 = 0$ and its root given by Equation (2.31) is considered for approximating value of $y = y_1 + a_1$, where $y_1 > 0$. Further correction is applied by way of $a_2$ which is given by Equation (2.32) and then $y$ is considered equal to $y_3 = y_2 + a_2$ and the process is repeated. In this way, finally $a_n$ is given by Equation (2.35) and then $y$ is approximated as $y_{n+1} = y_n + a_n$. To avoid repetition, the details are not given in this section but can be had by referring to section 2.1.3.

### 3.3.1. Assignment Of Value To $y_1$

Since $y > 0$, and approximation involves $\ln(\ln X/y_1)$, that implies $y_1/X < 1$. It is also explicit from the graph that value of $y$ approximates $X$ particularly, when $X$ is large. In such a situation, $y_1$ is assigned a value approximating $X$ such that $y_1/X < 1$. For approximation of value $y$, when $X > 15$, the calculator must have high precision as in such cases value of $y$ tends to approach value $X$ and assignment of value to $y_1$ should be close to $X$ so that $y_n + a_n$ should neither be a negative quantity nor more than $X$.

On the basis of the formulae derived in section 2.1.3, the approximation of value $y$ is given in Table 3. No error was found in approximated value of $y$ up to nine decimal points, when correction was applied four times or less.

Table 3 Displaying correctness of quadratic equation approximation of $y$ in $-ye^{e^{-y}} = -X$

| $X$ | $y_1$ | $y_2 = y_1 + a_1$ | $y_3 = y_2 + a_2$ | $y_4 = y_3 + a_3$ | $y_5 = y_4 + a_4$ | Actual val. |
|---|---|---|---|---|---|---|
| .001 | .0005 | $3.73099569 \times 10^{-4}$ | $3.68014666 \times 10^{-4}$ | $3.68014826 \times 10^{-4}$ | $3.680148264 \times 10^{-4}$ | $3.680148264 \times 10^{4}$ |
| 1 | 0.5 | 0.5676084521 | 0.5671432902 | 0.5671432904 | No need | 0.567143290 |
| 7 | 6.999 | 6.965264136 | 6.989105456 | 6.993661344 | 6.993578652 | 6.993578652 |
| 13 | 12.9999 | 12.99997596 | 12.99997059 | 12.99997061 | No need | 12.99997061 |

### 4. Playfulness of Equation $N(X) = y$ or $ye^{e^y} = X$ When Real Value Of $X$ Is Given

It is observed, while approximating value of $y$ satisfying equation $ye^{e^y} = X$, when given $X$ varies from zero to infinity that assignment of value from 1 to $10^{10}$ to $z_1$, where $z = \ln \ln y$, and $X = 10^3$ leads to the same precise result as $1.84021218$ irrespective of initial assumption of $z_1$.

It is also observed, while approximating value of $y$ satisfying equation $ye^{e^y} = X$, when given $X$ varies from zero to infinity that assignment of value from 1 to $10^8$ to $z_1$, where $z = \ln y$, and $X = 10^2$ leads to the same precise result as 1.44402854 irrespective of initial assumption of $z_1$. These characteristics of quadratic equation approximation methods engage one with curiosity how these methods lead to same destination of precise value of $y$. Tables 4 and 5, prove, starting with $z_1$ equal to 1 or $10^{10}$ in one case and $10^8$ in another, value of $y = \ln \ln (z_n + a_n)$ or $y = \ln (z_n + a_n)$ as the case may be, is same as 1.84021218 or 1.44402854. It is further submitted, the calculating device used for this paper is Desmos Scientific calculator and owing to its own limitation of precision, maximum value $z_1$ equal to $10^{10}$ in one case and $10^8$ in another case were obtained. Had the calculator been highly precise, value to $z_1$ can be assigned still higher. It is also important to note that value of $z_1$ does not matter much in these types of approximations, therefore, inadvertent human error that may creep in, will get self remedied. This self corrective action will be explained in detail in appropriate section.

It is also worth mentioning that assignment of value to $z_1$ does not matter match but it matters in approximation of $y$ in equation $-ye^{e^{-y}} = -X$. In such cases, there is not much margin for assumption of value of $y_1$ particularly when $X > 15$ and value of $y_1$ tends to value $X$. For high precision approximation, equally high precision calculator is needed.

Table 4 Displaying same result for $y = \ln \ln (z_{n+1})$ in equation $z \ln(\ln z) = 10^3$

| $X$ | $z_1$ | $z_2$ | $z_3$ | $z_4$ | $z_5$ | $z_6$ | $z_7$ | $y = \ln(\ln z_8)$ | Actual |
|---|---|---|---|---|---|---|---|---|---|
| $10^3$ | 2 | 1229.791516 | 541.5325181 | 543.4155971 | No need | No need | No need | 1.84021218 | 1.84021218 |
| $10^3$ | $10^2$ | 568.3000447 | 543.4152668 | 543.4155969 | 543.41559693 | No need | No need | 1.84021217 | 1.84021218 |
| $10^3$ | $10^7$ | 377.62 | 543.5919838 | 543.4155969 | 543.41559693 | No need | No need | 1.84021217 | 1.84021218 |
| $10^3$ | $10^{10}$ | 10 | 700.4612574 | 543.3556818 | 543.4155969 | No need | No need | 1.84021217 | 1.84021218 |

Table 5 Displaying same result for $y = \ln (z_{n+1})$ in equation $z + \ln(\ln z) = \ln (10^2)$

| $X$ | $z_1$ | $z_2 = z_1 + a_1$ | $z_3 = z_2 + a_2$ | $z_4 = z_3 + a_3$ | $z_5 = z_4 + a_4$ | $y = \ln z_5$ | Actual val. |
|---|---|---|---|---|---|---|---|
| $10^2$ | 2 | 4.283737107 | 4.237733285 | 4.237733378 | No need | 1.444028546 | 1.44402854 |
| $10^2$ | $10^2$ | 3.58462743 | 4.238098722 | 4.237733377 | 4.237733378 | 1.4440285 | 1.4440285 |
| $10^2$ | $10^8$ | 1.01 | 7.241328797 | 4.228166831 | 4.237733378 | 1.4440285 | 1.4440285 |

## 5. Formulae Of N Function And Its Relationship With Lambert W Function

Lambert W Function $ye^y = x$ can be written $e^{ye^y} = e^x$ by raising both sides to base $e$. Substituting $y$ with Lambert W Function $W(x)$ and rearranging

$$ye^{e^y} = ye^{x/y} = W(x)e^{x/W(x)} = X \qquad (5.1)$$

or

$$N(X) = W(x) = y, \qquad (5.2)$$

in N notation. Please distinguish between $X$ given by Equation (5.1) and $x$. Similarly, it can be proved,

$$N\!\left(-W(x)e^{W(x)/x}\right) = -W(x), \qquad (5.3)$$

or

$$N\!\left(W(x)e^{x/W(x)}\right) + N\!\left(-W(x)e^{W(x)/x}\right) = 0.$$

When $x = 1$, using Equation (5.3) and (5.2)

$$N(-1) = -W(1) = -\Omega, \qquad (5.4)$$
$$N(\Omega e^{1/\Omega}) = W(1) = \Omega \qquad (5.5)$$

where $\Omega$ is Omega constant. Taking $y = 1$, $y = -1$, $y = e$, $y = \pi i$, $y = -\pi i$, it can be proved,

$$N(e^e) = W(e) = 1, \qquad (5.6)$$
$$N(-e^{1/e}) = W(-1/e) = -1, \qquad (5.7)$$
$$N(e^{1+e^e}) = W(e^{e+1}) = e, \qquad (5.8)$$
$$N(\pi i/e) = W(-\pi i) = \pi i,$$

$$N(-\pi i/e) = W(\pi i) = -\pi i.$$

In general, if $y = b\,i$, then $bi\,e^{e^{bi}} = -be^{\cos(b)}\{\sin(\sin b) - i\cos(\sin b)\}$ or
$$N(X) = b\,i,$$
where $X = -be^{\cos(b)}\{\sin(\sin b) - i\cos(\sin b)\}$.

### 5.1. N Function Constant

N function constant is that value of $y$ which makes $N(1) = y$ or $ye^{e^y} = 1$. Value of $y$ satisfying $N(1) = y$ has already been calculated in Table 1 and 2. It is 0.2698741376 precise up to ten places of decimals. It is also that value of $y$ which satisfies the equation $ln(y) + e^y = 0$. Interestingly, referring to Equation (5.4), value of $y$ in equation $N(-1) = y$, equals minus Omega constant, mathematically, $N(-1) = -\Omega$.

## 6. Solution To General N Function $y^p e^{e^y} = X$

In this method, it is considered $e^y = z$ and equation $y^p e^{e^y} = X$, where $p$ is a real number positive or negative, transforms to $(ln\,z)^p e^z = X$, which can be written as
$$z + p\,ln\,ln\,z = ln\,X. \qquad (6.1)$$

Let $z = z_1 + a_1$, then Equation (6.1) takes the form
$$z_1 + a_1 + p\,ln\,ln\,(z_1 + a_1) = ln\,X. \qquad (6.2)$$

Following the procedure given in section 2.1.2, Equation (6.2) takes the form, $a_1^2 - l_1 a_1 - m = 0$, where
$$l_1 = -2(p + z_1 ln\,z_1)/(1 + ln\,z_1) + (ln\,X - p\,ln\,ln\,z_1 - z_1) \qquad (6.3)$$
$$m_1 = 2z_1 ln\,z_1 (ln\,X - p\,ln\,ln\,z_1 - z_1)/(1 + ln\,z_1). \qquad (6.4)$$

and the root $a_1 = (1/2)\left(l_1 + \sqrt{l_1^2 + 4m_1}\right)$ of the above said quadratic equation is considered for approximation when $p$ is positive real number. Therefore,
$$z = z_2 = z_1 + a_1 = (1/2)\left(l_1 + \sqrt{l_1^2 + 4m_1}\right),$$

The correction is iterated $n$ times and value of
$$z = z_1 + a_1 + a_2 + a_3 + \cdots + a_n$$
approximates $y$ by equation $y = ln\,z$, where
$$a_n = (1/2)\left(l_n + \sqrt{l_n^2 + 4m_n}\right) \qquad (6.5)$$
$$l_n = -2(p + z_n ln\,z_n)/(1 + ln\,z_n) + (ln\,X - p\,ln\,ln\,z_n - z_n), \qquad (6.6)$$
$$m_n = 2z_n ln\,z_n (ln\,X - p\,ln\,ln\,z_n - z_n)/(1 + ln\,z_n). \qquad (6.7)$$

With each iteration,
$$(ln\,X - p\,ln\,ln\,z_1 - z_1) > (ln\,X - p\,ln\,ln\,z_2 - z_2) > \cdots > (ln\,X - p\,ln\,ln\,z_n - z_n)$$
therefore, $|a_1| > |a_2| > |a_3| > \cdots > |a_n|$ and $z$ converges to $z_1 + a_1 + a_2 + a_3 + \cdots + a_n$.
To avoid repetition, details are not given in this section but can be accessed at sections 2.1.2 and 2.2.2. But when $p$ is negative real number, the. equation $y^p e^{e^y} = X$, takes the form $y^{-P} e^{e^y} = X$, where $p = -P$ and $P$ is a positive real number. Examination of this equation reveals, when value of $y$ increases, value of $y^{-P}$ decreases and value of $e^{e^y}$ increases so as to satisfy equation $y^{-P} e^{e^y} = X$. That means one lower and one higher value of $y$ satisfy the equation $y^{-P} e^{e^y} = X$. In other words equation $y^{-P} e^{e^y} = X$ has two solutions for y. For its other solution, roots
$$a_1' = (1/2)\left(l_1 - \sqrt{l_1^2 + 4m_1}\right)$$
$$a_2' = (1/2)\left(l_2' - \sqrt{l_2'^2 + 4m_2'}\right)$$
$$\ldots$$
$$a_n' = (1/2)\left(l_n' - \sqrt{l_n'^2 + 4m_n'}\right)$$

are considered and $l_n^{'}$, $m_n^{'}$ and $z_n^{'}$ are given by equations
$$l_n^{'} = -2(p + z_n^{'} \ln z_n^{'})/(1 + \ln z_n^{'}) + (\ln X - p \ln \ln z_n^{'} - z_n^{'})$$
$$m_n^{'} = 2z_n^{'} \ln z_n^{'}(\ln X - p \ln \ln z_n^{'} - z_n^{'})/(1 + \ln z_n^{'}),$$
$$z = z_{n+1} = z_1 + a_1^{'} + a_2^{'} + a_3^{'} + \cdots + a_n^{'}$$

To prove the veracity of the formulae derived, some examples are given in Table 6.

Table 6 Displaying veracity of quadratic equation approximation of z in $z + p \ln(\ln z) = \ln X$

| $X$ | $p$ | $z_1$ | $z_2 = z_1 + a_1$ | $z_3 = z_2 + a_2$ | $z_4 = z_3 + a_3$ | $z_5 = z_4 + a_4$ | $y = \ln z_5$ | Actual val. |
|---|---|---|---|---|---|---|---|---|
| $10^{-3}$ | 100 | 2 | 2.491102068 | 2.485184811 | 2.4851848164 | No need | 0.9103470296 | 0.9103470296 |
| 1 | 1/2 | 2 | 1.033635801 | 1.130170443 | 1.113808277 | 1.1138082775 | 0.1077850239 | 0.1077850239 |
| $10^5$ | 10 | 2 | 6.511463632 | 5.836013178 | 5.836419631 | No need | 1.764117532 | 1.764117532 |
| $10^5$ | $10^3$ | 2 | 2.762445334 | 2.742333232 | 2.742333407 | No need | 1.008809166 | 1.008809166 |
| $10^5$ | $10^{-5}$ | 2 | 11.51291905 | 11.51291655 | 11.51291654 | 11.51291653 | 2.443469582 | 2.443469582 |
| $10^{-6}$ | $-100$ | 2 | 3.371424416 | 3.275194137 | 3.275205452 | 3.2752054516 | 1.1863806 | 1.1863806 |
| $10^{-6}$ | $-100$ | 2 | 64.64745923 | 145.8953508 | 146.9267621 | 146.92676283 | 4.989934251 | 4.989934251 |
| 50 | $-5$ | 2 | 1.97027827 | 1.970268079 | 1.970280603 | 1.97028060384 | 0.6781759711 | 0.6781759711 |
| 50 | $-5$ | 2 | 6.377805661 | 7.370745345 | 7.371966944 | 7.3719669445 | 1.997684556 | 1.997684556 |
| 10 | $-1/100$ | 2 | 2.300787216 | 2.300760756 | 2.3007607554 | No need | 0.8332398315 | 0.8332398315 |
| 75 | $-1/5$ | 2 | 4.965580494 | 4.39601952 | 4.395990129 | 4.3959901296 | 1.480692791 | 1.480692791 |

When the equation is of the form $(-y)^p e^{e^{-y}} = -X$, where $y$ has positive real value, then the equation can be written $-y^p e^{e^{-y}} = -X$, when $p$ is of the form $(2n_1 + 1)$ or $(2n_1 + 1)/(2n_2 + 1)$ and $n_1$ and $n_2$ are real integers. But when $p$ is of the form $2n_1$ or $2n_1/(2n_2 + 1)$ or $(2n_1 + 1)/2n_2$, then equation $(-y)^p e^{e^{-y}} = -X$ does not have solution for $y$ in real domain. Equation $(-y)^p e^{e^{-y}} = -X$ or $y^p e^{e^{-y}} = X$ transforms to quadratic equation $a_1^2 - l_1 a_1 - m = 0$, where

$$l_1 = -\{y_1 + \ln n(X/y_1^p)\} + 2\{y_1 \ln(X/y_1^p) - p\}/\{1 - \ln X/y_1^p\}, \quad (6.8)$$
$$m_1 = 2y_1 \ln(X/y_1^p)\{y_1 + \ln \ln(X/y_1^p)\}/\{1 - \ln(X/y_1^p)\} \quad (6.9)$$

and $y_1$ is the assigned value such that $y_1 < X^{1/p}$. Using iterations as already explained, value of y can be approximated. Based on the above equations some examples are given Table 7.

Table 7 Displaying veracity of quadratic equation approximation of $y$ in $y_1 + \ln \ln(X/y_1^p) = 0$

| $X$ | $p$ | $y_1$ | $y_2 = y_1 + a_1$ | $y_3 = y_2 + a_2$ | $y_4 = y_3 + a_3$ | $y_5 = z_4 + a_4$ | $y_6 = z_5 + a_5$ | Actual val. |
|---|---|---|---|---|---|---|---|---|
| $10^{-10}$ | 10 | .09 | .09136201169 | .09127692626 | .09127653083 | 0.09127652716 | No need | 0.09127652716 |
| 1 | $-10$ | 1.01 | 1.022214158 | 1.032852467 | 1.035960542 | 1.03611954 | 1.036119908 | 1.0361199078 |

## 7. Applications

N Function $N(X) = y$, where value real value of $X$, is applicable to physical processes that involve equation $ye^{e^y} = X$ or other equations transformable to this equation. Some such equations are given hereinafter.

### 7.1 Solution To Equations Transformable To $ye^{e^y} = X$

Processes involving equations

$$p + \ln z + e^{qz} = r, \qquad (7.1)$$
$$z = r \ln \ln(pq/z) \qquad (7.2)$$
$$\ln \ln z + z = p \qquad (7.3)$$
$$z \ln \ln z = p \qquad (7.4)$$
$$e^z \ln z = p \qquad (7.5)$$
$$aye^{e^{by}} = p \qquad (7.6)$$

are transformable to equation $ye^{e^y} = X$ and hence solvable by the methods as explained in the paper. Equation (7.1) takes the form $ye^{e^y} = X$, when $y = (b\,z)$, $X = (be^{r-p})$. Equation (7.2) takes the form $ye^{e^y} = X$, when $y = (z/r)$ and $X = (pq/r)$. Equation (7.3) takes the form $ye^{e^y} = X$, when $y = \ln z$ and $X = e^p$. Equation (7.4) takes the form $ye^{e^y} = X$, when $y = \ln \ln z$ and $X = p$. Equation (7.5) takes the form $ye^{e^y} = X$, when $y = \ln z$ and $X = p$. Equation (7.6) takes the form $ye^{e^y} = X$, when $by = z$ and $X = bp/a$.

### 7.2. Solution To Equations Transformable To General N Function $y^p e^{e^y} = X$

Processes involving equations

$$q \ln z + e^{rz} = s, \tag{7.7}$$
$$z = r \ln \ln (p/z^q) \tag{7.8}$$
$$p \ln \ln z + z^q = r \tag{7.9}$$
$$p z + r e^{e^{z^q}} = s \tag{7.10}$$
$$e^{z^p} \ln qz = r \tag{7.11}$$
$$z \ln \ln (p/z^q) = r \tag{7.12}$$
$$z \, e^{e^{z^q}} = p \tag{7.13}$$

are transformable to equation $y^p e^{e^y} = X$ and hence solvable by the methods as explained in the paper. Equation (7.7) takes the form $y^q e^{e^y} = X$, when $y = (r\,z)$, $X = r^q e^s$. Equation (7.8) takes the form $y e^{e^y} = X$, when $y = (z/r)$ and $X = (pr^q)$. Equation (7.9) takes the form $y^p e^{e^y} = X$, when $y = q \ln z$ and $X = q^p e^r$. Equation (7.10) takes the form $(p/r) \ln \ln y + y^q = s/r$ when $z = \ln \ln y$. This equation is same as Equation (7.9) and therefore can be solved as already explained. Equation (7.11) takes the form $y^P/q^p + \ln \ln y = \ln r$, when $y = qz$ and this equation is same as Equation (7.9) and its solution has already been explained. Equation (7.12) takes the form $y^{-q} e^{e^y} = X$, when $r/z = y$ and $X = p/r^q$. Equation (7.13) takes the form $y^{1/q} e^{e^y} = p$, when $y = z^q$. It can thus be summed up that an equation that involves terms 1) $z$ and $e^{e^z}$ or 2) $\ln z$ and $e^z$ or 3) $\ln \ln z$ and $z$ in algebraic sum or product form can be solved by the methods explained in the paper.

## 8. Comparison Of This Method With Newton Method Of Approximation

Newton's method uses the iteration

$$x_{n+1} = x_n - f(x_n)/f'(x), \tag{8.1}$$

where $f(x_n)$ and $f'(x_n)$ are the values of the function and its derivative at $x = x_n$. Initially,

$$x_1 = x_0 - f(x_0)/f'(x_0)$$

and $x_0$ is assigned value. Care is taken while assigning value to $x_0$ so that it is the in the vicinity of the root of $f(x) = 0$ to yield a convergent series [8]. The method also requires that $f(x_n)$ is contiguous differentiable and does not equal zero. If $f(x)$ has roots with multiplicity more one then one, then the convergence rate is merely linear unless special care is taken [8]. The method devised in this paper uses iteration given by Equation (2.14) which is reproduced below

$$z_{n+1} = z_n + a_n = z_n + (1/2)\left(l_n + \sqrt{l_n^2 + 4m_n}\right)$$

Initially,

$$z_2 = z_1 + a_1 = z_1 + (1/2)\left(l_1 + \sqrt{l_1^2 + 4m_1}\right)$$

and $z_0$ is assigned value, $l_1$ and $m_1$ are coefficient of z and constant term of quadratic equation

$$a_1^2 - a_1 l_1 - m_1 = 0$$

obtained from equations $y e^{e^y} = X$, when $\ln z = e^y$ or $z = e^y$ as the case may be. These iterations based on quadratic approximation do not involve derivative of function $f(z)$, hence it is not the requirement that $f(z)$ should be continuously differentiable, it is not the requirement that $f'(z) \neq 0$. More importantly, unlike Newton's Method, it is not the requirement that assignment to value of $z_1$ should be carefully done so as to obtain a convergent series. Referring to section 3.1.1, assignment of value to $z_1$ can be made from 1 to $10^8$ when $z = \ln y$, and $y e^{e^y} = 10^2$ and it leads to the same precise result as 1.44402854 irrespective of initial assumption of $z_1$. Tables 4 and 5, prove, starting with $z_1$ equal to 1 or $10^{10}$ in one case and $10^8$ in another, value of $y = \ln \ln (z_n + a_n)$ or $y = \ln (z_n + a_n)$ as the case may be, is same as 1.84021218 or 1.44402854. However, both methods provide sequences $x_{n+1}$ and $a_{n+1}$ that converge super linearly to zero.

## 9. Results And Conclusions

*There are processes in physical and living sciences, which involve functions containing i) a variable and double exponential of the variable or ii) logarithm of a variable and exponential of the variable or iii) logarithm of logarithm of a variable and also the variable itself.* Denying existence of such processes would be more fallacy than truthfulness. Such processes can be characterised by an equation $ye^{e^y} = X$ or in more general form $y^p e^{e^y} = X$, where given $X$ and $p$ have real values except plus infinity, minus infinity or zero. Equation $ye^{e^y} = X$ is given a nomenclature of N Function and is written as $N(X) = y$, where value of $X$, is given and value of $y$ is approximated. *In a published paper [7], it is proved, $\ln(1 + 1/t)$ approximates $2/(2t + 1)$ provided $t \gg 1$.* If $t$ happens to be not far greater than 1, then the said approximation will not be precise and will need further correction which is applied by way of iteration. Aforesaid approximation relates to logarithm but I need approximation of $y$ in equations $ye^{e^y} = X$ and $-ye^{e^{-y}} = -X$, when given $X$ is a positive quantity. Equation $ye^{e^y} = X$ can be written in the following forms i and ii whereas $-ye^{e^{-y}} = -X$ in form iii.

i. $z \ln \ln z = X$ where $\ln z = e^y$,
ii. $z + \ln \ln z = \ln X$, where $z = e^y$,
iii. $y = -\ln \ln(X/y)$.

Assuming $z = z_1 + a_1$ or $y = y_1 + a_1$, where $z_1$ or $y_1$ (as the case may be,) are assigned the values arbitrarily and applying the approximation $\ln(1 + 1/t) \simeq 2/(2t + 1)$, the above said equations i, ii, iii yield quadratic equations in $a_1$. Roots of the quadratic equation facilitate determining the value of $z_1 + a_1$ or $y_1 + a_1$. Let $z_1 + a_1$ or $y_1 + a_1$ be denoted by $z_2$ or $y_2$ respectively. Since initial assumption of value of $z_1$ or $y_1$ may not be close to $z$ or $y$ (as the case may be), therefore, further corrections, by way of iterations assuming $z = (z_2 + a_2), z = (z_3 + a_3), \ldots, z = (z_n + a_n)$, are needed. While assuming value of $z_1$, care will have to be taken that $z_n + a_n$ or $y_n + a_n$ is a real positive quantity. If it is not so, then assumption of value of $z_1$ or $y_1$ will have to be changed. Value of integer $n$ should be such that the result is of desired precision. Higher the value of $n$, more precise will be the result. In this paper, maximum value of $n$ used is 6 and six times iterations yielded the result precise up to nine decimal places.

Invariably root $l_n/2 + (1/2)(l_n^2 + 4m_n)^{1/2}$ of quadratic equation $a_n^2 - a_n l_n - m_n = 0$ gave positive values of $z_n + a_n$ or $y_n + a_n$ and was thus used and other discarded. But in equation $y^p e^{e^y} = X$, for solution of $y$, when $p$ was a negative quantities, both roots $l_n/2 + (1/2)(l_n^2 + 4m_n)^{1/2}$ and $l_n/2 - (1/2)(l_n^2 + 4m_n)^{1/2}$ yielded positive values of $z_n + a_n$ or $y_n + a_n$, hence both were thus used.

Some interesting relationships between Lambert W Function $W(x) = y$ or $ye^y = x$ and N Function $N(X) = y$ or $ye^{e^p} = X$, were noticed. When $X = W(x)e^{x/W(x)}$, then $N(X) = W(x)$. Also $N(-W(x)e^{W(x)/x}) = -W(x)$, When $x = 1$, then $-W(1)e^{W(1)} = -1$ and $N(-1) = -W(1) = -\Omega$, where $\Omega$ is Omega constant. For the sake of establishing relationship between the two functions, capital $X$ was used for equation $ye^{e^y} = X$. Value of N Function constant i.e. $N(1)$ was approximated 0.2698741376, which is precise to nine decimal places.

General N Function $y^p e^{e^y} = X$ and $-y^p e^{e^{-y}} = X$, where p is a real quantity except plus infinity, minus infinity or zero, can be written in the following form iv and v respectively.

iv. $z + p \ln \ln z = \ln X$, where $y = \ln$,
v. $y_1 + \ln \ln (X/y_1^p) = 0$.

These equations on application of approximation $\ln(1 + 1/t) \simeq 2/(2t + 1)$ and writing $z = z_1 + a_1$ or $y = y_1 + a_1$, transforms to quadratic equations. Roots of these equation facilitate approximation of z or y. Iterating the correction $n$ times as explained in equation $ye^{e^y} = X$, precise value of $z$ or $y$ can be approximated and from value of $z$, value of $y$ is calculated as $ln\ z$.

### 8.1. Precision

In equation $-ye^{e^{-y}} = -X$, for solution of $y$ particularly, when $X$ is large, value of $y$ approximates $X$. But value of $y$ must be slightly less than $X$ on account of fact value of $e^{e^{-y}}$ is more than one. Therefore, for precise results, it is required, calculating device must be of high precision. For example, Desmos Scientific Calculator being used for computation purposes for this paper, yields result $y = X$, when given $X \geq 23$, whereas precise result of $y$ should be less than $X$.

### 8.2. Self Corrective Nature

It is explicit from this method, value of $a_1$ is approximated on the basis of assumption of value of $z_1$ or $y_1$. Value of $a_2$ is approximated on the basis of assumption of value of $z_1 + a_1$ say $z_2$ or $y_1 + a_1$ say $y_2$. If there is an inadvertent error committed by the operator in recording $z_1$ or $a_1$, then value of $a_2$ will be approximated on the basis $z_1 + a_1$ i.e. taking into consideration inadvertent error. At last, approximation of $a_n$ will be based upon $z_{n-1} + a_{n-1}$ or $z_1 + a_1 + a_2 + a_3 + \cdots + a_{n-1}$ i.e. taking into consideration all earlier correction, therefore, $a_n$ will eliminate earlier inadvertent errors to yield $z = z_{n+1} = z_n + a_n$. Thus the method has inherent self correction and in doing so, the number of total correction i.e, number of iterations may increase.

### 8.3. Recreational Nature

For positive $X$ in equation $ye^{e^y} = X$ written as $z\ ln(ln\ z) = X$, when $y = ln(ln\ z)$, or as $z + \ln(ln\ z) = \ln X$. when $y = \ln z$, assuming different and varying values of $z_1$ yield precise result of $y$. When $X = 10^3$ in equation $z\ ln(ln\ z) = X$, assumption of value of $z_1$ from 2 to $10^{10}$ yielded precise value of 1.84021218 in each case. When $X = 10^2$ in equation $z + \ln(ln\ z) = \ln X$, assumption of value of $z_1$ from 2 to $10^8$ yielded precise value of 1.4440285 in each case.

**Disclosure of Conflict of Interest, Statements and Declarations Competing Interests, Funding and Acknowledgement**

*Disclosure Conflict of Interest, Competing interest, Financial or Non-financial intestinal,*

There is no competing interest, there is no financial or non financial interest involved in this research work. There is also no conflict of interest with anyone. However, I acknowledge the help provided by the website https://www.desmos.com in calculating the values of tedious exponential and logarithmic terms.

*Acknowledgement of Funds, if any*

No fund, whatsoever in cash, kind or any financial support has been received for this work. However, I acknowledge the help provided by the website https://www.desmos.com in calculating the values of tedious exponential and logarithmic terms.

*Ethical statement*

This statement is not applicable to this research work.

*Consent statement*

This statement is not applicable to this research work.

**Author's Contribution**

Entire research paper is the sole contribution of the corresponding author.